\documentclass[12pt]{amsart}

\usepackage[dvips]{graphicx}
	\usepackage{amsmath,amssymb,amsthm,wrapfig,amsfonts,enumerate,latexsym}

	\newtheorem{dfn}{Definition}[section]
	\newtheorem{thm}[dfn]{Theorem}
	\newtheorem{prop}[dfn]{Proposition}
	\newtheorem{cor}[dfn]{Corollary}
	\newtheorem{lem}[dfn]{Lemma}
	\newtheorem{rem}[dfn]{Remark}
	
 	\newtheorem{claim}[dfn]{Claim}

	\newtheorem{ack}{Acknowledgements\!\!}

	\newcounter{yon}
	\numberwithin{equation}{section}

	\def\notin{\not\in}

	\newcommand{\dist}{\mathop{\mathit{d}} \nolimits}
	
	\newcommand{\diam}{\mathop{\mathrm{diam}} \nolimits}

	\newcommand{\sep}{\mathop{\mathrm{Sep}} \nolimits}
	
	\newcommand{\obs}{\mathop{\mathrm{ObsDiam}}  \nolimits}

	\newcommand{\obsc}{\mathop{\mathrm{ObsCRad}}               \nolimits}
	\newcommand{\crad}{\mathop{\mathrm{CRad}}           \nolimits}

	\newcommand{\obin}{\mathop{\mathrm{Obs}L^p\mathrm{\text{-}Var}}
	\nolimits}
	\newcommand{\obinin}{\mathop{\mathrm{Obs}L^2\mathrm{\text{-}Var}}
	\nolimits}

\begin{document}

	\title[Central and $L^p$-concentration of $1$-Lipschitz maps into $\mathbb{R}$-trees]
	{Central and $L^p$-concentration of $1$-Lipschitz maps\\ into
	$\mathbb{R}$-trees}
	\author[Kei Funano]{Kei Funano}
	\address{Mathematical Institute, Tohoku University, Sendai 980-8578, JAPAN}
	\email{sa4m23@math.tohoku.ac.jp}
	\subjclass[2000]{53C21, 53C23}
	\keywords{median, mm-space, observable $L^p$-variation, observable diameter, observable
	central radius, $\mathbb{R}$-tree}
	\thanks{This work was partially supported by Research Fellowships of
	the Japan Society for the Promotion of Science for Young Scientists.}
	\dedicatory{}
	\date{\today}

	\maketitle


\begin{abstract}In this paper, we study the L\'{e}vy-Milman
 concentration phenomenon of $1$-Lipschitz maps from mm-spaces to
 $\mathbb{R}$-trees. Our main theorems assert that the concentration to
 $\mathbb{R}$-trees is equivalent to the concentration to the real line. 
 \end{abstract}
	\setlength{\baselineskip}{5mm}

    \section{Introduction}
    This paper is devoted to investigating the L\'{e}vy-Milman concentration
    phenomenon of $1$-Lipschitz maps from mm-spaces (metric measure spaces) to
    $\mathbb{R}$-trees. Here, an mm-space is a triple
    $(X,\dist_X,\mu_X)$ of a set $X$, a complete separable
    distance function $\dist_X$ on $X$, and a finite Borel measure
    $\mu_X$ on $(X,\dist_X)$. Let
    $\{(X_n,\dist_{X_n},\mu_{X_n})\}_{n=1}^{\infty}$ be a sequence of mm-spaces and
    $\{(Y_n,\dist_{Y_n})\}_{n=1}^{\infty}$ a sequence of metric spaces. Given a sequence $\{f_n
    :X_n \to Y_n\}_{n=1}^{\infty}$ of $1$-Lipschitz maps, we consider
    the following three properties:

    (i) (Concentration property) There exist points $m_{f_n}\in Y_n$, $n\in \mathbb{N}$, such that
     \begin{align*}
      \mu_{X_n}\big(\{  x_n \in X_n \mid \dist_{Y_n}(f_n(x_n),m_{f_n})
      \geq \varepsilon          \}\big) \to 0 \text{ as }n\to \infty
      \end{align*}for any $\varepsilon >0$.

      (ii) (Central concentration property) The maps $f_n$, $n\in \mathbb{N}$, concentrate to the
      center of mass of the push-forward measure $(f_n)_{\ast}(\mu_{X_n})$. In other words, the concentration
      property (i) holds in the case where $m_{f_n}$ is the center of
      mass.

       (iii) ($L^p$-concentration property) For a number $p>0$, we have 
       \begin{align*}
        \int\int_{X_n \times X_n} \dist_{Y_n} \big( f_n
        (x_n),f_n(y_n)\big)^p \ d\mu_{X_n}(x_n ) d\mu_{X_n}(y_n) \to 0
        \text{ as }n\to \infty 
        \end{align*}Each target metric space $Y_n$, $n\in \mathbb{N}$,
        is called a \emph{screen}. Chebyshev's inequality proves
        that the $L^p$-concentration (iii)
        implies the concentration property (i) for any $p>0$. If each screen $Y_n$, $n\in
        \mathbb{N}$, is an Euclidean space $\mathbb{R}^k$, then the $L^p$-concentration
        (iii) for $p\geq 1$ yields
        the central concentration property (ii) (see Lemma \ref{l8}). The
        central concentration (ii) is stronger than the concentration
        property (i). There is an example of maps $f_n$, $n\in \mathbb{N}$, with
        the concentration property (i), but not having the central
        concentration property (ii) (see Remark \ref{remla}). In some
        special cases, the concentration (i) implies the central and
        $L^p$-concentration properties (ii) and (iii) (see \cite[Subsection 2.4]{funano1}
        and \cite[Section $3\frac{1}{2}.31$]{gromov}). 

        Vitali D. Milman first introduced the concentration and the central concentration properties (i) and (ii) for $1$-Lipschitz functions (i.e.,
        $Y_n=\mathbb{R}$, $n\in \mathbb{N}$) and emphasized their
        importance in his investigation of asymptotic geometric
        analysis (see \cite{mil}). Nowadays those properties are widely studied by many
        literature and blend with various areas of mathematics (see
        \cite{gromov}, \cite{ledoux}, \cite{mil2},
	\cite{milsch}, \cite{sch}, \cite{tal}, \cite{tal2} and references therein for
	further information). M. Gromov first considered the case of general
        screens in 
        \cite{gromovcat}, \cite{gromov2}, and \cite[Chapter
        $3\frac{1}{2}$]{gromov}. See \cite{funano1}, \cite{funano}, and \cite{ledole} for another works of general
        screens. In
        \cite{gromov}, Gromov settled the concentration and central
        concentration properties (i) and
        (ii) for $1$-Lipschitz maps by
        introducing the
        \emph{observable diameter} $\obs_Y(X;-\kappa) $ and the
        \emph{observable central radius} $\obsc_Y(X;-\kappa)$ for an
        mm-space $X$, a metric space $Y$, and $\kappa >0$ (see Section
        $2$ for the
        precise definitions). The $L^2$-concentration property (iii) was
        first appeared in Gromov's paper \cite{gromovcat}. Motivated by \cite{gromovcat}, the author
        introduced in \cite{funano1} the \emph{observable
        $L^p$-variation} $\obin_Y(X)$ to study the property (iii) (see
        Section 2 for the definition). Note that given a sequence $\{
        X_n\}_{n=1}^{\infty}$ of mm-spaces and $\{ Y_n\}_{n=1}^{\infty}$
        of metric spaces, $\obs_{Y_n}(X_n;-\kappa) $ (resp.,
        $\obsc_{Y_n}(X_n;-\kappa)$, $\obin_{Y_n}(X_n)$) converges to zero as $n\to \infty$ for any $\kappa >0$ if and
        only if any sequence $\{ f_n:X_n \to Y_n \}_{n=1}^{\infty}$ of $1$-Lipschitz
        maps (resp., central, $L^p$-)concentrates.


In this paper, we treat the case of $\mathbb{R}$-tree screens. 
\begin{thm}\label{p1}Let $\{  X_n
 \}_{n=1}^{\infty}$ be a sequence of mm-spaces. Then, the following
 (\ref{lala1}) and (\ref{lala2}) are equivalent to each other.
 \begin{align}\label{lala1}
  \obs_{\mathbb{R}} (X_n;-\kappa)\to 0 \text{ as }n\to \infty \text{ for
  any }\kappa >0.
  \end{align} 
	 \begin{align}\label{lala2}
	  \sup \{\obs_T(X_n;-\kappa)    \mid T     \text{ is
	  an }\mathbb{R}\text{-tree}               \}\to 0 \text{ as }n\to
	  \infty \text{ for any }\kappa >0.
	  \end{align}
     \end{thm} 
Theorem \ref{p1} is a complete solution
      to Gromov's exercise in \cite[Section
      $3\frac{1}{2}.32$]{gromov}. In \cite[Section
    5]{funano1}, the author proved it only for
    simplicial tree screens. The implication (\ref{lala2}) $\Rightarrow$
      (\ref{lala1}) is obvious. For the proof of the converse, we define the notion of
      a \emph{median} for a finite Borel measure on an $\mathbb{R}$-tree in
      Section $3$ and
      proves that any $1$-Lipschitz maps $f_n$ from $X_n$ into
      $\mathbb{R}$-trees concentrate to medians for the push-forward
      measure $(f_{n})_{\ast}(\mu_{X_n})$.

      To study the central and $L^p$-concentration for (ii) and (iii) into
      $\mathbb{R}$-trees, we estimate the distance between the center of
      mass and a median of a finite Borel measure on an $\mathbb{R}$-tree from the
      above in Section 5. For this estimate, we partially extend K-T. Sturm's characterization of the center of mass
      on a simplicial tree to the case of an $\mathbb{R}$-tree (see Proposition
      \ref{p7} and Section 4). From the estimate, we bound $\obsc_{T}(X;-\kappa)$ (resp., $\obin_{T}(X)$)
      from the above in terms of $\obsc_{\mathbb{R}}(X;-\kappa)$ (resp.,
      $\obin_{\mathbb{R}}(X)$) (see Propositions \ref{c3} and
      \ref{c4}). As a result, we obtain
        \begin{thm}\label{t1}Let $\{ X_n\}_{n=1}^{\infty}$ be a sequence 
         of mm-spaces. Then, the following (\ref{lala3}) and
         (\ref{lala5}) are equivalent to each other.
 \begin{align}\label{lala3}
  \obsc_{\mathbb{R}}(X_n;-\kappa) \to 0 \text{ as }n\to \infty \text{
  for any }\kappa >0.
  \end{align}
 \begin{align}\label{lala5}
  \sup \{ \obsc_{T}(X_n;-\kappa) \mid T \text{ is an }\mathbb{R}\text{-tree}
  \} \to 0\text{ as }n\to \infty \text{ for any }\kappa >0.
  \end{align}
      \end{thm}

      \begin{thm}\label{t2}Let $\{ X_n\}_{n=1}^{\infty}$ be a sequence
       of mm-spaces and $p\geq 1$. Then, the following (\ref{lala4}) and (\ref{lala6}) are equivalent to
       each other.
 \begin{align}\label{lala4}
  \obin_{\mathbb{R}}(X_n) \to 0 \text{ as }n\to \infty.
  \end{align}
 \begin{align}\label{lala6}
  \sup \{ \obin_{T}(X_n) \mid T \text{ is an }\mathbb{R}\text{-tree}
  \} \to 0\text{ as }n\to \infty.
  \end{align}
       \end{thm}

       The condition (\ref{lala3}) is stronger than (\ref{lala1}) (see
       Lemma \ref{l7} and Remark \ref{remla}), and (\ref{lala4}) implies
       (\ref{lala3}) (see Lemma \ref{l8}). It seems that the conditions
       (\ref{lala3}) and (\ref{lala4}) are not equivalent, but we have no counterexample. 

      In our previous work, the author investigated the above properties
        (i), (ii), and (iii) for $1$-Lipschitz maps into Hadamard
        manifolds (see \cite[Theorems 1.3, 1.4, and Lemma 4.4]{funano1}). The $L^2$-concentration property
        (iii) in that case
        is also studied by Gromov (see \cite[Section 13]{gromovcat}). Our theorems are thought as of $1$-dimensional analogue to these works. 
     \section{Preliminaries}
     \subsection{Basics of the concentration and the $L^p$-concentration}
     \subsubsection{Observable diameter and separation distance}
     Let $Y$ be a metric space and $\nu$ a Borel measure on $Y$ such that $m:=\nu(Y)<+\infty$. 
	 We define for any $\kappa >0$
	 \begin{align*}
	  \diam (\nu , m-\kappa):= \inf \{ \diam Y_0 \mid Y_0 \subseteq Y \text{ is a Borel subset such that }\nu(Y_0)\geq m-\kappa\}
	  \end{align*}and call it the \emph{partial diameter} of $\nu$.

     \begin{dfn}[Observable diameter]\upshape Let $(X,\dist_X,\mu_X)$ be
      an mm-space with $m:=\mu_X(X)$ and $Y$ a metric space. For any $\kappa >0$ we
	 define the \emph{observable diameter} of $X$ by 
	 \begin{align*}
	  \obs_Y (X; -\kappa):=
	   \sup \{ \diam (f_{\ast}(\mu_X),m-\kappa) \mid f:X\to Y \text{ is a
      }1 \text{{\rm -Lipschitz map}}  \}. 
      \end{align*}The target metric space $Y$ is called the \emph{screen}.
      \end{dfn}
      The idea of the observable diameter comes from the quantum and statistical
	mechanics, that is, we think of $\mu_X$ as a state on a configuration
	space $X$ and $f$ is interpreted as an observable. 

     Let $(X,\dist_X,\mu_X)$ be an
      mm-space. For any $\kappa_1,  \kappa_2\geq 0$, we
      define the \emph{separation distance} $\sep
      (X;\kappa_1,\kappa_2)= \sep (\mu_X;\kappa_1,\kappa_2)$ of $X$ as the supremum of
      the distance $\dist_X(A,B)$, where $A$ and $B$ are Borel subsets
      of $X$ satisfying that $\mu_X(A)\geq \kappa_1$ and $\mu_X(B)\geq \kappa_2$.

      The proof of the following lemmas are easy and we omit the proof.
      \begin{lem}[{cf.~\cite[Section $3\frac{1}{2}.33$]{gromov}}]\label{l1}Let
       $(X,\dist_X,\mu_X)$ and $(Y,\dist_Y,\mu_Y)$ be two
       mm-spaces. Assume that a $1$-Lipschitz map $f:X\to \mathbb{R}$
       satisfies $f_{\ast}(\mu_X)=\mu_Y$. Then we have
       \begin{align*}
        \sep (Y;\kappa_1,\kappa_2)\leq \sep (X;\kappa_1,\kappa_2)
        \end{align*}
       \end{lem}

       \begin{lem}\label{l2}For any $\kappa >m/2$, we have $\sep (X;\kappa,\kappa)=0$.
        \end{lem}

       The relationships between the observable diameter and the
       separation distance are the following:

       \begin{prop}[cf.~{\cite[Section $3\frac{1}{2}.33$]{gromov}}]\label{p2}Let $(X,\dist,\mu)$ be an mm-space and $0<\kappa'< \kappa$. Then we have
	 \begin{align*}
	  \sep (X;\kappa,\kappa)\leq \obs_{\mathbb{R}}(X;-\kappa') .
     \end{align*}
        \end{prop}

	\begin{prop}[cf.~{\cite[Section $3\frac{1}{2}.33$]{gromov}}]\label{p3}
	 For any $\kappa >0$, we have
	 \begin{align*}
	   \obs_{\mathbb{R}}(X;-2\kappa) \leq \sep (X;\kappa,\kappa).
	  \end{align*}
        \end{prop}
        See \cite[Subsection $2.2$]{funano} for details of the proofs of the above propositions.

       \begin{cor}[{cf.~\cite[Section $3\frac{1}{2}.33$]{gromov}}]\label{c1}A
        sequence $\{ X_n\}_{n=1}^{\infty}$ of mm-spaces satisfies that
        \begin{align*}
         \obs_{\mathbb{R}}(X_n;-\kappa)\to 0 \text{ as }n\to \infty
         \end{align*}for any $\kappa>0$ if and only if $\sep (X_n;\kappa ,\kappa)\to 0$ as $n\to
        \infty$ for any $\kappa >0$.
        \end{cor}

        \subsubsection{Observable $L^p$-variation}

         Let $(X,\dist_X,\mu_X)$ be an mm-space and $(Y,\dist_Y)$ a metric
space. Given a Borel measure $\nu$ on $Y$ and $p\in (0,+\infty)$, we put
\begin{align*}
 V_p(\nu):=  \Big(\int \int_{Y\times Y}
 \dist_Y(x,y)^p   \ d\nu(x) d \nu(y)\Big)^{1/p}.
 \end{align*}For a Borel measurable map $f:X\to Y$, we also put
 $V_p(f):= V_p \big(f_{\ast}(\mu_X)\big)$

Let $ \{  X_n        \}_{n=1}^{\infty}$ be a
 sequence of mm-spaces and $\{  Y_n   \}_{n=1}^{\infty}$ a sequence of
 metric spaces. For any $p\in (0,+\infty]$, we say that a sequence $\{ f_n :X_n \to Y_n
 \}_{n=1}^{\infty}$ of Borel measurable maps \emph{$L^p$-concentrates}
 if $V_{p}(f_n)\to 0$ as $n\to \infty$.

 Given an mm-space $X$ and a metric space $Y$ we define
\begin{align*}
\obin_Y(X):= \sup \{   V_p(f) \mid f:X\to Y \text{ is a }1
 \text{-Lipschitz map}                  \},
\end{align*}and call it the \emph{observable} $L^p$\emph{-variation} of
 $X$.

 \begin{lem}\label{l3}For any closed subset $A\subset X$, we have
  \begin{align*}
   \obin_{\mathbb{R}}(A) \leq \obin_{\mathbb{R}}(X).
   \end{align*}
  \begin{proof}Let $f:A\to \mathbb{R}$ be an arbitrary $1$-Lipschitz
   function.  From \cite[Theorem 3.1.2]{ambro}, there exists a $1$-Lipschitz extension of
 $f$, say $\widetilde{f}:X\to \mathbb{R}$. Hence, we get
   \begin{align*}
    V_p(f)\leq V_p (\widetilde{f}) \leq \obin_{\mathbb{R}}(X).
    \end{align*}This completes the proof.
   \end{proof}
  \end{lem}

See \cite[Subsection 2.4]{funano1} for the relationships between the
observable diameter and the observable $L^p$-variation.
        \subsection{Basics of $\mathbb{R}$-trees}
        Before reviewing the definition of $\mathbb{R}$-trees, we recall
    some standard terminologies in metric geometry. Let $(X,\dist_X)$ be
    a metric space. A rectifiable curve $\eta:[0,1]\to X$ is called a
    \emph{geodesic} if its arclength coincides with
    the distance $\dist_X\big(\eta(0),\eta(1)\big)$ and it has a constant speed,
    i.e., parameterized proportionally to the arc length. We say that
    $(X,\dist_X)$ is a \emph{geodesic space} if any two points in $X$ are joined
    by a geodesic between them. Let $X$ be a geodesic space. A \emph{geodesic
    triangle} in $X$ is the union of the image of three geodesics joining a triple of points in $X$ pairwise. A subset $A\subseteq X$
    is called \emph{convex} if every geodesic joining two points in $A$
    is contained in $A$.

    A complete metric space $(T,\dist_T)$ is called an 
    \emph{$\mathbb{R}$-tree} if it has the following properties:
    \begin{itemize}
		   \item[$(1)$]For all $z,w\in T$ there exists a unique unit
                       speed geodesic
                       $\phi_{z,w}$ from $z$ to $w$.
		   \item[$(2)$]The image of every simple path in $T$ is the
                       image of a geodesic.
	 \end{itemize}Denote by $[z,w]_T$ the image of the
     geodesic $\phi_{z,w}$. We also put $(z,w]_T:=[z,w]_T\setminus
     \{z\}$ and $(z,w)_T:=[z,w]_T\setminus \{z,w\}$. A complete geodesic space $T$
     is an $\mathbb{R}$-tree if and only if it is $0$-hyperbolic,
     that is to say, every edge in any geodesic triangle in $T$ is included in
     the union of the other two edges. See \cite{chiswell} for another characterizations of $\mathbb{R}$-trees.
     Given $z\in T$, we indicate by $\mathcal{C}_T(z)$ the set of all connected
    components of $T\setminus \{ z\}$. We also denote by $\mathcal{C}_T'(z)$
     the set of all $\{z\}\cup T'$ for $T'\in \mathcal{C}_T(z)$. Although the following lemma is
     somewhat standard, we prove it for the completeness. 
     \begin{lem}\label{l4}Each $T'\in \mathcal{C}_T(z)$ is convex.
  \begin{proof}From the property $(2)$
   of $\mathbb{R}$-trees, it is sufficient to prove that $T'$ is arcwise
   connected. Taking a point $z\in T'$, we put
   \begin{align*}A:= \{ w\in T' \mid z
   \text{ and }w \text{ are connected by a path in }T'\}.
    \end{align*}It is easy to see that the set $A$ is closed in
   $T'$. Since every metric ball in $T$ is arcwise connected, the set $A$
   is also open. Since $T'$ is connected, we get $T'=A$. This completes the proof.
  \end{proof}
     \end{lem}

     A subset in an $\mathbb{R}$-tree is called
    a \emph{subtree} if it is a closed convex subset. Note that a
    subtree is itself an $\mathbb{R}$-tree.



    \begin{prop}\label{p4}Every connected subset in an $\mathbb{R}$-tree is convex.
     \begin{proof}Let $T$ be an $\mathbb{R}$-tree. Suppose that there exists a connected subset
      $T'\subseteq T$ which is not convex. Then, there are
      points $z,w\in T'$ and $\widetilde{z}\in (z,w)_T$ such
      that $\widetilde{z}\notin T'$. Since $T' =
      \bigcup \{T'\cap C \mid C\in \mathcal{C}_T(\widetilde{z})\}$ and each
      $C\in \mathcal{C}_T(\widetilde{z})$ is open, from the
      connectivity of $T'$, there is $C_0\in \mathcal{C}_T(\widetilde{z})$ such that $T'\subseteq
      C_{0}$. Since $C_{0}$ is convex by Lemma \ref{l4},
      we get $\widetilde{z}\in [z,w]_T \subseteq
      C_{0}$. This is a contadiction since $\widetilde{z}\notin
      C_{0}$. This completes the proof.
      \end{proof}
    \end{prop}

\subsection{Center of mass of a measure on a CAT(0)-space and
  observable central radius}

  \subsubsection{Basics of the center of mass of a measure on CAT(0)-spaces}
In this subsection, we review Sturm's works about measures on a
CAT(0)-spaces. Refer \cite{jost} and \cite{sturm} for details. A geodesic metric
space $X$ is called a CAT$(0)$-space if we have
\begin{align*}
\dist_X\big(x,\gamma (1/2)\big)^2 \leq \frac{1}{2}\dist_X(x,y)^2 +
 \frac{1}{2}\dist_X(x,z)^2 - \frac{1}{4} \dist_X(y,z)^2
\end{align*}for any $x,y,z \in X$ and any minimizing geodesic $\gamma
 :[0,1]\to X$ from $y$ to $z$. For example, Hadamard manifolds, Hilbert
 spaces, and $\mathbb{R}$-trees are all CAT(0)-spaces.

Let $(X,\dist_X)$ be a metric space. We denote by $\mathcal{B}(X)$ the set of all
finite Borel measures $\nu$ on $X$ with the separable support. We
indicate by $\mathcal{B}^1(X)$ the set of all Borel measures $\nu \in \mathcal{B}(X)$ such
that $\int_{X}\dist_X(x,y) \ d\nu (y)<+\infty$ for
some (hence all) $x\in X$. We also indicate by $\mathcal{P}^1(X)$ the
set of all probability measures in $\mathcal{B}^1(X)$. For any $\nu \in \mathcal{B}^1(X)$ and $z\in X$, we consider the function
$h_{z,\nu}:X\to \mathbb{R}$ defined by
\begin{align*}
h_{z,\nu}(x):= \int_{X} \{\dist_X(x,y)^2-\dist_X(z,y)^2 \} \ d\nu(y).
\end{align*}
Note that
\begin{align*}
\int_X |\dist_X(x,y)^2-\dist_X(z,y)^2| \ d\nu(y)\leq \dist_X(x,z)\int_X
 \{\dist_X(x,y)+ \dist_X(z,y)\} \ d \nu(y) <+\infty.
\end{align*}A point $z_0 \in X$ is called the
 \emph{center of mass} of the measure $\nu \in \mathcal{B}^1(X)$ if for
 any $z\in X$, $z_0$ is a unique minimizing
 point of the function $h_{z,\nu}$. We denote the point $z_0$ by $c(\nu)$.
A metric space $X$ is said to be \emph{centric} if every $\nu \in
 \mathcal{B}^1(X)$ has the center of mass. 

\begin{prop}[{cf.~\cite[Proposition
 $4.3$]{sturm}}]\label{p6}
A CAT(0)-space is centric.
\end{prop}

A simple variational argument yields the following lemma.
\begin{lem}[{cf.~\cite[Propsition $5.4$]{sturm}}]\label{l5}
Let $H$ be a Hilbert space. Then for each $\nu \in \mathcal{B}^1(H)$ with $m=\nu(X)$, we have
\begin{align*}
c(\nu)= \frac{1}{m}\int_H y \ d\nu(y).
\end{align*}
\end{lem}
Let $(T,\dist_T)$ be an $\mathbb{R}$-tree and $\nu \in \mathcal{B}^1(T)$. For $z\in T$ and $T'\in
\mathcal{C}_T'(z)$, we put
\begin{align*}
 c_{z,T'}(\nu):= \int_{T'}\dist_T (z,w) \ d\nu (w) - \int_{T\setminus
 T'}\dist_T (z,w) \ d\nu(w).
 \end{align*}

 Let us consider a (possibly infinite) simplicial tree $T_s$. Here, the
 length of each
 edge of $T_s$ is not necessarily equal to $1$. We assume that every
 vertex of $T_s$ is an isolated point in the vertex set of $T_s$. 
\begin{prop}[{cf.~\cite[Proposition $5.9$]{sturm}}]\label{p7}
 Let $\nu \in \mathcal{B}^1(T_s)$ and $z\in T_s$. Then, $z=c(\nu)$ if and only if $c_{z,T'}(\nu)\leq 0$ for any $T'\in \mathcal{C}_{T_s}'(z)$.
 \end{prop}

 \begin{prop}[{cf.~\cite[Proposition $6.1$]{sturm}}]\label{p8}Let $N$ be a CAT(0)-space and $\nu \in \mathcal{B}^1(N)$. Assume that
  the support of $\nu$ is contained in a closed convex subset $K$ of
  $N$. Then, we have $c(\nu)\in K$. 
  \end{prop}

  Let $X$ be a metric space. For $\mu,\nu \in \mathcal{P}^1(X)$, we
  define the \emph{$L^1$-Wasserstein distance} $\dist^W_1(\mu,\nu)$ between $\mu$
  and $\nu$ as the infimum of $\int_{X\times X}\dist_X(x,y) \ d\pi
  (x,y)$, where $\pi \in \mathcal{P}^1(X\times X)$ runs over all \emph{couplings}
  of $\mu$ and $\nu$, that is, the measures $\pi$ with the propery that $\pi (A\times
  X)=\mu(A)$ and $\pi (X\times A)=\nu (A)$ for any Borel subset
  $A\subseteq X$. 
  \begin{lem}[{cf.~\cite[Theorem 7.12]{villani}}]\label{l6}A sequence $\{ \mu_n \}_{n=1}^{\infty} \subseteq
   \mathcal{P}^1(X)$ converges to $\mu \in \mathcal{P}^1(X)$ with
   respect to the distance function $\dist_1^W$ if and only if the sequence $\{
   \mu_n\}_{n=1}^{\infty}$ converges weakly to the measure $\mu$ and
   \begin{align*}
    \lim_{n\to \infty}\int_X \dist_X (x,y) \ d\mu_n(y) = \int_X \dist_X
    (x,y) \ d\mu (y)
    \end{align*}for some (and then any) $x\in X$.
   \end{lem}

  \begin{thm}[{cf.~\cite[Theorem $6.3$]{sturm}}]\label{t4}Let $N$ be a
   CAT(0)-space. Given $\mu,\nu \in \mathcal{P}^1(N)$, we have
    $\dist_N(c(\mu),c(\nu))\leq \dist^{W}_1(\mu, \nu)$.
   \end{thm}

    \subsubsection{Observable central radius}

    Let $Y$ be a metric space and assume that $\nu \in \mathcal{B}^1(Y)$ has the
center of mass. We denote by $B_Y(y,r)$ the closed ball in $Y$ centered at
$y\in Y$ and with raidus $r>0$. For any $\kappa >0$, putting $m:= \nu (Y)$, we define the \emph{central radius}
$\crad(\nu,m-\kappa)$ of $\nu$ as the infimum of $\rho >0$ such that
$\nu\big(B_Y(c(\nu),\rho)\big)\geq m-\kappa$.

Let $(X,\dist_X, \mu_X)$ be an mm-space with $\mu_X \in \mathcal{B}^1(X)$ and $Y$
a centric metric space. For any $\kappa >0$, we define 
\begin{align*}
\obsc_Y(X;-\kappa):= \sup \{ \crad(f_{\ast}(\mu_X),m-\kappa) \mid f:X\to
 Y \text{ is a }1 \text{-Lipschitz map}\},
\end{align*}and call it the \emph{observable central radius} of $X$.

\begin{lem}[{cf.~\cite[Section $3\frac{1}{2}.31$]{gromov}}]\label{l7}For any $\kappa >0$, we have
\begin{align*}
\diam (\nu,m-\kappa)\leq 2\crad(\nu,m-\kappa).
\end{align*}In particular, we get
 \begin{align*}
  \obs_Y (X;-\kappa)\leq 2\obsc_Y(X;-\kappa).
  \end{align*}
\end{lem}

\begin{rem}\label{remla}\upshape From the above lemma, we see that the
 central concentration implies
 the concentration. The converse is not true in general. For
 example, consider the mm-spaces $X_n:= \{x_n,y_n \}$ with distance
 function $\dist_{X_n}$ given by $\dist_{X_n}(x_n,y_n):=n$ and with a
 Borel probability measure $\mu_{X_n}$ given by $\mu_{X_n}(\{
 x_n\}):=1-1/n$ and $\mu_{X_n}(\{ y_n\}):=1/n$. Then, $1$-Lipschitz maps
 $f_n:X_n \to \mathbb{R}$ defined by $f_n(x):= \dist_{X_n}(x,x_n)$
 satisfy that
 \begin{align*}
  (f_n)_{\ast}(\mu_{X_n})\big(B_{\mathbb{R}}(c((f_n)_{\ast}(\mu_{X_n})),1/2)\big)=0
  \end{align*}for any $n\in \mathbb{N}$, whereas
 $\obs_{\mathbb{R}}(X_n;-\kappa) \to 0$ as $n\to \infty$.
 \end{rem}

\begin{lem}\label{l8}Let $\nu \in
 \mathcal{B}^1(\mathbb{R}^n)$ with $m:=\nu(\mathbb{R}^n)$. Then, for any
 $p\geq 1$ and $\kappa >0$, we have
 \begin{align}\label{exc1}
  \crad(\nu,m-\kappa) \leq \frac{V_p (\nu)}{(m\kappa)^{1/p}}.
  \end{align}In the case of $p=2$, we also have the better estimate
 \begin{align}\label{exc2}
  \crad(\nu,m-\kappa)\leq \frac{V_2(\nu)}{\sqrt{2m\kappa}}.
  \end{align}
 \begin{proof}We shall prove that $\nu \big(\mathbb{R}^n\setminus
 B_{\mathbb{R}^n}\big(c(\nu),\rho_0      \big)\big)\leq \kappa$ for $\rho_0 :=
 V_p(\nu)/(m\kappa)^{1/p}$. Suppose that $\nu \big(\mathbb{R}^n\setminus
 B_{\mathbb{R}^n}\big(c(\nu),\rho_0      \big)\big)>\kappa$. From Lemma \ref{l5}, we get
  \begin{align*}
   \int_{\mathbb{R}^n}|c(\nu)-x|^p \ d\nu (x) \leq \frac{V_p(\nu)^p}{m}. 
  \end{align*}Hence, from Chebyshev's inequality,
 we see that
\begin{align*}
\frac{V_p(\nu)^p}{m}=\rho_0^p \kappa <\int_{\mathbb{R}^n} |c(\nu)-x|^p \
 d \nu (x) \leq \frac{V_p(\nu)^p}{m},
\end{align*}which is a contradiction. Therefore, we obtain $\nu \big(
 B_{\mathbb{R}^n}(c(\nu)  ,\rho_0       )\big)\geq m-\kappa$ and so (\ref{exc1}).

  Since
  \begin{align*}
   \int_{\mathbb{R}^n}|c(\nu)-x|^2 \ d\nu(x)=\frac{V_2(\nu)^2}{2m},
   \end{align*}the same argument yields (\ref{exc2}).
  This completes the proof.
\end{proof}
 \end{lem}

 \begin{cor}Let $X$ be an mm-space with $\mu_X \in
 \mathcal{B}^1(X)$. Then, for any $p\geq 1$, we have
  \begin{align}\label{exc3}
   \obsc_{\mathbb{R}^n}(X;-\kappa) \leq \frac{1}{(m\kappa)^{1/p}}\obin_{\mathbb{R}^n}(X).
   \end{align}In the case of $p=2$, we also have the better estimate
 \begin{align}\label{exc4}
  \obsc_{\mathbb{R}^n}(X;-\kappa)\leq  \frac{1}{\sqrt{2m\kappa}}\obinin_{\mathbb{R}^n}(X)
 \end{align}
  \end{cor}
 \begin{cor}\label{c2}Let $X$ be an mm-space. Then, for any $p\geq 1$
  and $\kappa >0$, we have
  \begin{align}\label{exc5}
   \sep (X;\kappa,\kappa) \leq \frac{2}{(m\kappa)^{1/p}}\obin_{\mathbb{R}}(X). 
   \end{align}
  In the case of $p=2$, we also have
  \begin{align}\label{exc6}
   \sep (X;\kappa ,\kappa)\leq \sqrt{\frac{2}{m\kappa}} \obinin_{\mathbb{R}}(X).
   \end{align}
  \begin{proof}Assume first that there is a $1$-Lipschitz function $f:X\to \mathbb{R}$
   such that $f_{\ast}(\mu_X) \notin
   \mathcal{B}^1(\mathbb{R})$. From H\"{o}lder's inequality, we have
   $\int_{\mathbb{R}} |x-y|^p \ df_{\ast}(\mu_X)(y) =+\infty$ for any
   $x\in X$. This implies $V_p(f)=+\infty$ and so  $\obin_{\mathbb{R}}(X)=+\infty$.

   We consider the other case that $f_{\ast}(\mu_X)\in
   \mathcal{B}^1(\mathbb{R})$ for any $1$-Lipschitz function $f:X\to
   \mathbb{R}$. Combining Proposition \ref{p2} with Lemma \ref{l7} and (\ref{exc3}), we have
   \begin{align*}
    \sep (X;\kappa,\kappa) \leq \frac{2}{(m\kappa')^{1/p}} \obin_{\mathbb{R}}(X)
    \end{align*}for any $\kappa > \kappa'>0$. Letting $\kappa'\to
   \kappa$, we have (\ref{exc5}). Replacing (\ref{exc3}) with (\ref{exc4}) in the above argument, we also obtain (\ref{exc6}). 
   \end{proof}
  \end{cor}

  \section{Existence of a median on an $\mathbb{R}$-tree}
Let $T$ be an $\mathbb{R}$-tree and $\nu$ a finite Borel measure on $T$
     with $m:=\nu(T)<+\infty$. A median of $\nu$ is a point
     $z\in T$ such that there exist two subtrees $T',T''\subseteq T$ such
     that $T=T'\cup T''$, $T'\cap T''=\{ z\}$, $\nu (T')\geq m/3$, and
     $\nu(T'')\geq m/3$. The existence of a median of a finite Borel measure on a simplicial tree is proved in
     \cite[Proposition 5.2]{funano1}. The purpose of this section is to prove the existence of a median of a
  finite Borel measure on an $\mathbb{R}$-tree, which is needed for the
  proofs of our main theorems. Although the proof of the existence is similar to the proof for the case of a simplicial tree, we prove it for the completeness: 
      \begin{prop}\label{pp1}Every finite Borel measure on an $\mathbb{R}$-tree has
      a median.
      \begin{proof}Let $\nu$ be a finite Borel measure on an
       $\mathbb{R}$-tree with $m:=\nu(T)$. Assume that a point $z\in T$ satisfies that
       $\nu({T}')<m/3$ for any ${T}'\in C_T'(z)$, then
       it is easy to check that $z$ is a median of $\nu$. So, we assume that for any $z\in T$ there exists $T(z)\in
       \mathcal{C}_T'(z)$ such that $\nu (T(z))\geq m/3$. If for some $z\in T$,
       there exists $T'\in \mathcal{C}_T'(z)\setminus  \{T(z) \}$ such that $\nu
       (T')\geq m/3$, then this $z$ is a median of $\nu$. Thereby, we
       also assume that $\nu (T')<m/3$ for any $z\in T$ and $T'\in
       \mathcal{C}_T'(z)\setminus \{ T(z)\}$.

       Fixing a point $z_0\in T$, we assume that there exists
       $z \in T({z_0})\setminus \{z_0\}$ such that $z_0 \in T(z)$. Put
       \begin{align*}t_0 := \inf
       \{ t\in (0,\dist_T(z_0,z)]  \mid z_0 \in
        T(\phi_{z_0,z}(t))\}.
        \end{align*}
       \begin{claim}\label{claim1}$\phi_{z_0,z}(t_0)$ is a median of $\nu$.
        \begin{proof}Assume first that $t_0 =0$. Then, taking a
         monotone decreasing sequence $\{ t_n\}_{n=1}^{\infty}\subseteq (0,\dist_T(z_0,z)]$ such
         that $t_n \to 0$ as $n\to \infty$ and $z_0 \in
         T(\phi_{z_0,z}(t_n))$ for any $n\in \mathbb{N}$, we shall show that $\bigcap_{n=1}^{\infty}T(\phi_{z_0,z}(t_n))\subseteq
         \big(T\setminus T(z_0)\big)\cup \{z_0\}$. If it is, we conclude
         that the point $z_0
         = \phi_{z_0,z}(0)$ is a median of $\nu$ as follows: From the
         uniqueness of $T(\phi_{z_0,z_1}(t_n))$, we have $T(\phi_{z_0,z}(t_{n+1})) \subseteq
         T(\phi_{z_0,z}(t_n))$ for each $n\in \mathbb{N}$. Thus,
         we get
          $\nu \big(
          \bigcap_{n=1}^{\infty}T(\phi_{z_0,z}(t_n))\big) =
          \lim_{n\to \infty} \nu (T(
          \phi_{z_0,z}(t_n)))\geq {m}/{3}$.

         Suppose that
         there exists $w\in T(z_0)\setminus \{ z_0\}\cap
         \bigcap_{n=1}^{\infty}T(\phi_{z_0,z}(t_n))$. Note
         that $(z_0,z]_T\cap (z_0,w]_T\neq \emptyset$. Actually,
         suppose that $(z_0,z]_T\cap (z_0,w]_T= \emptyset$. Then, it
         follows from the property $(2)$ of $\mathbb{R}$-trees that
         $[z,w]_T = [z_0,z]_T \cup [z_0,w]_T$. Especially, we have $z_0
         \in [z,w]_T$. Since $T(z_0)\setminus
         \{z_0\}$ is convex by virtue of Lemma \ref{l4}, $[z,w]_T$
         does not contain the point $z_0$. This is a contradiction. Thus,
         there exists $t\in (0,\dist_T(z_0,z)]$ such that
         $\phi_{z_0,z}(t)\in (z_0,z]_T\cap (z_0,w]_T$.
         We pick $n_0\in \mathbb{N}$ with $t_{n_0} <t$. Since $w\in
         T(z_0)\setminus \{z_0\}\cap \bigcap_{n=1}^{\infty}T(\phi_{z_0,z}(t_n)) \subseteq
         T(\phi_{z_0,z}(t_{n_0}))\setminus \{ z_0\}$, we get $\phi_{z_0,z}(t)\in
         (z_0,w]_T\subseteq T(\phi_{z_0,z}(t_{n_0}))\setminus
         \{z_0\}$. Thereby, we get $\phi_{z_0,z}(t)\in
         T(\phi_{z_0,z}(t_{n_0}))\setminus \{ \phi_{z_0,z}(t_{n_0})
         \}$. Therefore, since $z_0 \in
         T(\phi_{z_0,z}(t_{n_0}))\setminus \{ \phi_{z_0,z}(t_{n_0})\}$
         and $T(\phi_{z_0,z}(t_{n_0}))\setminus \{
         \phi_{z_0,z}(t_{n_0})\}$ is convex, we obtain
         \begin{align*}
          \phi_{z_0,z}(t_n)\in [z_0,\phi_{z_0,z}(t)]_T\subseteq
          T(\phi_{z_0,z}(t_n)) \setminus \{ \phi_{z_0,z} (t_n) \}.
          \end{align*}This is a contradiction. Therefore, we have $\bigcap_{n=1}^{\infty}T(\phi_{z_0,z}(t_n))\subseteq
         \big(T\setminus T(z_0)\big)\cup
         \{z_0\}$. 

         We consider the other case that $t_0 >0$. Take a monotone
         increasing sequence $\{ t_n\}_{n=1}^{\infty}\subseteq
         (0,+\infty)$ such that $t_n\to t_0$ as $n\to \infty$ and $z_0
         \notin T\big(\phi_{z_0,z}(t_n)\big)$ for each $n\in
         \mathbb{N}$. Then, the same proof in the case of $t_0=0$ implies
         that $\nu
         \big(\bigcap_{n=1}^{\infty}T(\phi_{z_0,z}(t_n))\big)\geq
         m/3$ and
         $\bigcap_{n=1}^{\infty}T(\phi_{z_0,z}(t_n))\subseteq
         \big( T\setminus T(\phi_{z_0,z}(t_0)) \big) \cup
         \{\phi_{z_0,z}(t_0)\}$. Therefore, $\phi_{z_0,z}(t_0)$ is a
         median of $\nu$. This completes the proof of the claim.
         \end{proof}
        \end{claim}
       We next assume that $z_0 \notin T(z)$
       for any $z\in T(z_0)$. We denote by $\Gamma$ the set of all
       unit speed geodesics $\gamma:[0,L(\gamma)]\to T({z_0})$ such that
       $\gamma (0)=z_0$ and $\gamma ([t,L(\gamma)]) \subseteq T({\gamma
       (t)} )$ for any $t\in [0,L(\gamma)]$. Because of the assumption, we easily see 

       \begin{claim}\label{claim2}For any $z\in T(z_0)$, we have $\phi_{z_0,z}\in
        \Gamma$.
        \end{claim}

       \begin{claim}\label{claim3}For any $\gamma ,\gamma' \in \Gamma$ with
        $L(\gamma)\leq L(\gamma')$, we have
        \begin{align*}
        [\gamma (0),\gamma(L(\gamma))]_T\subseteq
         [\gamma'(0),\gamma'(L(\gamma'))]_T.
         \end{align*}
        \begin{proof}Suppose that
         \begin{align*}t_0:=\sup \{  t\in [0,L(\gamma)]
         \mid [\gamma (0),\gamma(t)]_T\subseteq
         [\gamma'(0),\gamma'(L(\gamma'))]_T      \}<
         L(\gamma).
          \end{align*}Then, we have $\gamma (t)\notin
         [\gamma'(0),\gamma'(L(\gamma'))]_T   $ for any
         $t>t'$. Actually, if $\gamma (t)\in
         [\gamma'(0),\gamma'(L(\gamma'))]_T $, then we have
         $\gamma (t)=\gamma'(t)$. Thus, $[\gamma (t_0),\gamma (t)]_T =
         [\gamma'(t_0),\gamma'(t)]_T$ by the property $(2)$ of the
         $\mathbb{R}$-trees. Thereby, we get $[\gamma (0),\gamma (t)]_T \subseteq
         [\gamma'(0),\gamma'(L(\gamma'))]_T$. Since $t>t_0$, this
         contradicts the definition of $t_0$. Therefore, from the
         property $(2)$ of $\mathbb{R}$-trees, we have
         \begin{align}\label{sl1}
          [\gamma(L(\gamma)),\gamma'(L(\gamma))]_T
          = [\gamma (t_0),\gamma (L(\gamma))]_T \cup
          [  \gamma' (t_0), \gamma'(L(\gamma))]_T .
          \end{align}
         Since $\gamma, \gamma' \in \Gamma$, we have
         $\gamma(L(\gamma )), \gamma'(L(\gamma)) \in
         T(\gamma (t_0))\setminus \{ \gamma (t_0)\}$. So, from
         the convexity of $T(\gamma (t_0))\setminus \{ \gamma
         (t_0)\}$, we
         get $[\gamma (L(\gamma)),
         \gamma'(L(\gamma'))]_T \subseteq T(\gamma
         (t_0))\setminus \{ \gamma (t_0)\}$. This is a contradition,
         because $\gamma(t_0) \in [\gamma (L(\gamma)),
         \gamma'(L(\gamma'))]_T$ from $($\ref{sl1}$)$. This
         completes the proof of the claim.
         \end{proof}
        \end{claim}
       Putting $\alpha : = \sup \{ L(\gamma) \mid \gamma \in \Gamma\}$, we
       shall show that $\alpha <+\infty$. If $\alpha <+\infty$, we
       finish the proof of the proposition as follows: From the
       completness of $\mathbb{R}$-trees and Claim \ref{claim3}, there exists a unique $\gamma \in
       \Gamma$ with $L(\gamma)=\alpha$. We also note that $\alpha >0$ by
       Claim \ref{claim2}. Thus, there exists a monotone
       increasing sequence $\{ t_n\}_{n=1}^{\infty}$ of positive numbers
       such that $t_n \to \alpha$ as $n\to \infty$. We easily see that
       $T(\gamma(t_{n+1}))\subseteq  T(\gamma(t_n))$ for any $n\in
       \mathbb{N}$ and $\bigcap_{n=1}^{\infty} T(\gamma(t_n))= \{ \gamma
       (L(\gamma)) \}$. Since $\nu \big(T(\gamma(t_n)) \big)\geq m/3$,
       the point $\gamma (L(\gamma))$ is a median of $\nu$.

       Suppose that $\alpha
       =+\infty$. Then, taking a sequence $\{
       \gamma_n\}_{n=1}^{\infty}\subseteq \Gamma$ such that
       $L(\gamma_n)< L(\gamma_{n+1})$ for any $n\in \mathbb{N}$ and
       $L(\gamma_n)\to +\infty$ as $n\to \infty$, we obtain $\bigcap_{n=1}^{\infty} 
        T( \gamma_n (L(\gamma_n))) =\emptyset$. Since $T(\gamma_n(L(\gamma_n)))
       \subseteq T(\gamma_{n+1}(L(\gamma_{n+1})))$ for
       any $n\in \mathbb{N}$, we have
       \begin{align*}
        0=\nu \Big( \bigcap_{n=1}^{\infty} T
        (\gamma_n (L(\gamma_n)))\Big)=\lim_{n\to \infty}\nu \big(  T
        (\gamma_n (L(\gamma_n))) \big) \geq \frac{m}{3},
        \end{align*}which is a contradiction. This completes the proof of the proposition.
       \end{proof}
      \end{prop}

  \section{The necessity of Proposition \ref{p7} for $\mathbb{R}$-trees}
   In order to prove the main theorems, we extend the necessity of Proposition
   \ref{p7} for $\mathbb{R}$-trees:
   \begin{prop}\label{p9}Let $T$ be an $\mathbb{R}$-tree and $\nu \in
    \mathcal{B}^1(T)$. Then, we
    have $c_{c(\nu),T'}(\nu)\leq 0$ for any $T'\in
    \mathcal{C}_T'(c(\nu))$.
    \begin{proof}For simplicities, we assume that $\nu(T)=1$. We shall
     approximate the measure $\nu$ by a measure whose support lies on
     a simplicial tree in $T$. Given
     $n\in \mathbb{N}$, there exists a compact subset $K_n \subseteq T$
     such that $\nu (T\setminus K_n) < 1/n \text{ and } \int_{T\setminus K_n}
      \dist_T (c(\nu),w) \ d\nu (w)<1/n$. Take a $(1/n)$-net $\{
     z_i^{n}\}_{i=1}^{l_n}$ of $K_n$ with mutually different elements
     such that $\dist_T(c(\nu),z_1^n)<1/n$. We then take a sequence
     $\{ A_i^n\}_{i=1}^{l_n}$ of mutually disjoint Borel subset of $K_n$ such that
     $z_i^n\in A_i^n$, $\diam A_i^n \leq 1/n$, and $K_n =
     \bigcup_{i=1}^{l_n}A_i^n$. Define the Borel probability measure
     $\nu_n$ on $\{ z_i^n\}_{i=1}^{l_n}$ by $\nu (\{ z_1^n\}):=
     \nu(A_1^n)+ \nu (T\setminus K_n)$ and $\nu(\{ z_i^n\}):=
     \nu(A_i^n) $ for $i\geq 2$.
     \begin{claim}\label{cl1}$\dist_1^{W}(\nu_n,\nu)\to 0$ as $n\to \infty$.
      \begin{proof}We shall show that 
       \begin{align}\label{s5}
        \lim_{n\to \infty} \int_T \dist_T (c(\nu),w)\ d\nu_n(w) = \int_T
        \dist_T(c(\nu),w)\ d\nu(w).
        \end{align}Since
       \begin{align*}
        \int_T \dist_T(c(\nu),w) \ d\nu_n (w)=
        \sum_{i=1}^{l_n}\dist_T(c(\nu),z_i^n)\nu (A_i^n) +
        \dist_T(c(\nu),z_1^n)\nu (T\setminus K_n),
        \end{align*}we have
       \begin{align}\label{s6}
        \Big| \int_T \dist_T (c(\nu),w) \ d\nu_n(w)- \sum_{i=1}^{l_n}
        \dist_T(c(\nu),z_i^n)\nu (A_i^n) \Big|<\frac{1}{n}.
        \end{align}From $\diam A_i^n <1/n$, we get
       \begin{align}\label{s7}
        & \Big| \sum_{i=1}^{l_n}\dist_T(c(\nu),z_i^n)\nu (A_i^n) -  \int_{K_n}
        \dist_T (c(\nu),w) \ d\nu(w)                      \Big|  \\ =\ & \Big|
        \sum_{i=1}^{l_n} \int_{A_i^n}\big\{
        \dist_T(c(\nu),w)-\dist_T(c(\nu),z_i^n)\big\} \ d\nu(w)\Big| 
        \leq  \sum_{i=1}^{l_n} \int_{A_i^n} \dist_T(w,z_i^n) \ d\nu
        (w) < \frac{1}{n}. \tag*{} 
        \end{align}Hence, combining (\ref{s6}) with (\ref{s7}) and
       \begin{align*}
        \Big| \int_{K_n}\dist_T (c(\nu),w) \ d\nu(w) - \int_{T}\dist_T (c(\nu),w)
        \ d\nu(w)\Big| \leq \int_{T\setminus K_n} \dist_T (c(\nu),w) \
        d\nu(w)< \frac{1}{n},
        \end{align*}we obtain (\ref{s5}). The same way of the above proof shows that the sequence
       $\{ \nu_n \}_{n=1}^{\infty}$ converges weakly to the measure
       $\nu$. Therefore, by using Lemma \ref{l6}, this completes the proof of
       the claim.
       \end{proof}
      \end{claim}Applying Claim \ref{cl1} to Theorem \ref{t4}, we get $c(\nu_n) \to c(\nu)$
     as $n\to \infty$. Since the convex hull in $T$ of the set
     $\{z_i^n\}_{i=1}^{l_n}$ is a simplicial tree with finite vertex set
     and $c(\nu_n)$ is contained in the convex hull by Proposition \ref{p8}, it
     follows from Proposition \ref{p7} that
     $c_{\widetilde{T},c(\nu_n)}(\nu_n)\leq 0$ for any $\widetilde{T}\in
     \mathcal{C}_T'(c(\nu_n))$. Let $T'\in
     \mathcal{C}_T'(c(\nu))$. 

     Assume first that $c(\nu_n)\in T\setminus T'$ for infinitely many
     $n\in \mathbb{N}$. Then, taking $T_n \in
     \mathcal{C}_T'(c(\nu_n))$ with $T'\subseteq T_n$, we have
     \begin{align*}
     c_{T',c(\nu)}(\nu_n)\leq
      c_{T_n,c(\nu_n)}(\nu_n)+\dist_T(c(\nu_n),c(\nu))\leq
      \dist_T ( c(\nu_n),c(\nu)).
      \end{align*}Therefore, we
     obtain $c_{T',c(\nu)}(\nu)= \lim_{n\to \infty}c_{T',c(\nu)}(\nu_n)\leq 0$.

     We consider the other case that $c(\nu_n )\in T'$ for any $n\in \mathbb{N}$. Let $z_n \in [c(\nu),c(\nu_1)]_T$ be the
     unique point such that
     \begin{align*}
      \dist_T(z_n, c(\nu_n) )= \inf \{ \dist_T(z,
     c(\nu_n))   \mid z\in [c(\nu),c(\nu_1)]_T\}.
      \end{align*}By taking a subsequence, we may assume that
     $\dist_T(c(\nu), z_{n+1}) \leq \dist_T(
     c(\nu),z_n)$ for any $n\in \mathbb{N}$. For each $n\geq 2$, we
     take $T_n \in \mathcal{C}_T'(z_n)$ and $\widetilde{T}_n \in
     \mathcal{C}_T'(c(\nu_n))$ such that $c(\nu_1) \in T_n$ and
     $c(\nu_1)\in \widetilde{T}_n$. Observe that $T_n \subseteq
     T_{n+1}$. Since $T_n \subseteq
     \widetilde{T}_n$, we have
     \begin{align}\label{slala}c_{T_n, z_n} (\nu_n) \leq c_{\widetilde{T}_n, c(\nu_n)}(\nu_n)+
      \dist_{T}(c(\nu_n),z_n) \leq \dist_T(c(\nu_n),z_n).
      \end{align}We also easily see
     \begin{claim}\label{clala}$T' \setminus \{ c(\nu)\}= \bigcup_{n=2}^{\infty}
      T_n$.
      \end{claim}The same proof of Claim \ref{cl1} implies that
     \begin{align*}
      \sup \Big\{   \Big|\int_{A} \dist_T(z_n,w) d\nu_n(w)- \int_A\dist_T(z_n,
      w) d\nu(w)\Big| \mid A\subseteq T \text{ is a Borel subset}\Big\}
      \to 0 \text{ as }n\to \infty.
      \end{align*}
     Combining this with (\ref{slala}) and Claim \ref{clala}, we obtain
     \begin{align*}c_{T',c(\nu)}(\nu)=
      \lim_{n\to \infty} c_{T_n,z_n} (\nu)= 
     \lim_{n\to \infty} c_{T_n, z_n} (\nu_n) \leq \lim_{n\to
      \infty}\dist_{T}(c(\nu_n),z_n)  = 0.
      \end{align*}
     This completes
     the proof of the proposition.
    \end{proof}
    \end{prop}

    The author does not know whether the converse of Proposition
    \ref{p9} holds or not.

\section{Proof of the main theorems}

  Combining Proposition \ref{pp1} with the same proof of \cite[Lemma
  5.3]{funano1} implies the following proposition:
     \begin{prop}\label{p5}Let $T$ be an $\mathbb{R}$-tree and $\nu$ a
      finite Borel measure. Then, for any $\kappa
      >0$, we have
      \begin{align}\label{s4}
       \nu \Big( B_T\Big(m_{\nu} ,
       \sep\Big(\nu;\frac{m}{3},\frac{\kappa}{2} \Big) \Big)
       \Big)\geq m-\kappa,
       \end{align}where $m_{\nu}$ is a median of the measure
      $\nu$. In particular, letting $X$ be an mm-space, we have
      \begin{align}\label{staba1}
       \obs_T (X;-\kappa)\leq 2\sep \Big(X;\frac{m}{3},\frac{\kappa}{2}\Big).
       \end{align} 
      \end{prop}
      Proposition \ref{p5} together with Corollary \ref{c1} yields
      Theorem \ref{p1}. The following way to prove Theorem \ref{p1} is much easier and more straightforward than the above way, that is,
      to prove the existence of a median of a measure on
$\mathbb{R}$-trees.

\begin{proof}[Proof of Theorem \ref{p1}]Our goal is to prove the following inequality:
      \begin{align}\label{s16}
       \obs_T (X;-\kappa)\leq 2\sep\Big(X;\frac{\kappa}{3},\frac{\kappa}{3}\Big)+4\obs_{\mathbb{R}}(X;-\kappa)
      \end{align} for any $\kappa>0$. Let $f:X\to T$ be an arbitrary $1$-Lipschitz map. Fixing a
 point $z_0\in T$, we shall consider the function $g:T\to \mathbb{R}$
 defined by $g(z):=\dist_T(z,z_0)$. Since $g\circ f:X\to \mathbb{R}$ is
 the $1$-Lipschitz function, from the definition of the observable
 diameter, there is an interval $A=[s,t]\subseteq [0,+\infty)$ such
 that $\diam A\leq \obs_T (X;-\kappa)$ and $(g\circ
 f)_{\ast}(\mu_X)(A)\geq m-\kappa$. Observe that the set $g^{-1}(A)$ is the annulus $\{ z\in T\mid s\leq \dist_T(z,z_0)\leq t\}$. We denote by
 $\mathcal{C}$ the set of all connected components of the set $g^{-1}(A)\setminus
 \{z_0\}$. 
 \begin{claim}\label{cl2}Assume that $s>0$. Then, for any $T'\in \mathcal{C}$, we have $\diam T' \leq 2\diam A$.
  \begin{proof}Given any $z_1,z_2\in T'$, we shall show that
   $\phi_{z_0,z_1}(s)=\phi_{z_0,z_1}(s)$. Suppose that
   $\phi_{z_0,z_1}(s)\neq \phi_{z_0,z_1}(s)$. Then, putting $s_0:= \sup
   \{ t\in [0,+\infty) \mid \phi_{z_0,z_1}(t) = \phi_{z_0,z_1}(t)\}$, we
   have $s_0< s$. From the definition of $s_0$ and the property $(2)$ of
   $\mathbb{R}$-trees, we have $(\phi_{z_0,z_1}(s_0),z_1]_T\cap
   (\phi_{z_0,z_2}(s_0),z_2]_T =\emptyset$. Therefore, from the property
   $(2)$ of $\mathbb{R}$-trees, we get
   \begin{align*}
    [z_1,z_2]_T= [\phi_{z_0,z_1}(s_0),z_1]_T\cup [\phi_{z_0,z_1}(s_0),z_2]_T.
    \end{align*}Hence, since $T'$ is convex by virtue of Proposition \ref{p4}, the points $z_1$ and $z_2$ must be included in
   different components in $\mathcal{C}_T(\phi_{z_0,z_1}(s_0))$. This is
   a contradiction, since
   $T'= \bigcup\{C \cap T'\mid C\in \mathcal{C}_T(\phi_{z_0,z_1}(s_0))\}$ and $T'$
   is connected. Thus, we have
   $\phi_{z_0,z_1}(s)=\phi_{z_0,z_2}(s)$. Consequently, we obtain
   \begin{align*}
    \dist_T(z_1,z_2)\leq \dist_T (z_1,\phi_{z_0,z_1}(s)) +
    \dist_T ( \phi_{z_0,z_2}(s),z_2)\leq 2(t-s)\leq 2\obs_{\mathbb{R}}(X;-\kappa).
    \end{align*}This completes the proof of the claim.
  \end{proof}
 \end{claim}
 Assume first that $s\leq \sep (X;\kappa/3,\kappa/3)/2$. Since every path
 connecting two components in $\mathcal{C}$ must cross the point $z_0$, by Claim
 \ref{cl2}, we have
 \begin{align*}\diam (f_{\ast}(\mu_X),m-\kappa)\leq \diam g^{-1}(A)\leq \sep \Big(X;\frac{\kappa}{3},\frac{\kappa}{3}\Big) +
 4\obs_{\mathbb{R}}(X ;-\kappa).
  \end{align*}

 We consider the other case that $s> \sep (X;\kappa /3, \kappa /3) /
 2$. Suppose that $f_{\ast}(\mu_X)(T')<\kappa /3$ for any $T'\in
 \mathcal{C}$. Since $f_{\ast}(\mu_X)(g^{-1}(A))\geq m-\kappa\geq \kappa$,
 we have $\mathcal{C}'\subseteq \mathcal{C}$ such that
 \begin{align*}
  \frac{\kappa}{3} \leq f_{\ast}(\mu_X)\Big(\bigcup \mathcal{C}'\Big)<\frac{2\kappa}{3}.
  \end{align*}Hence, by putting $\mathcal{C}'':=\mathcal{C}\setminus \mathcal{C}'$, we get
 \begin{align*}
  \sep \Big(X;\frac{\kappa}{3},\frac{\kappa}{3}\Big) <
  \dist_T\Big(\bigcup \mathcal{C}', \bigcup \mathcal{C}''\Big) \leq \sep
  \Big(f_{\ast}(\mu_X);\frac{\kappa}{3},\frac{\kappa}{3}\Big) \leq \sep \Big(X;\frac{\kappa}{3},\frac{\kappa}{3}\Big),
  \end{align*}which is a contradiction. Thereby, there exists $T'\in \mathcal{C}$
 such that $f_{\ast}(\mu_X)(T')\geq \kappa/3$. For a subset $A\subseteq
 T$ and $r>0$, we put $A_r:= \{ z \in T \mid \dist_T(z,A)\leq r \}$.
 \begin{claim}\label{cl3}$f_{\ast}(\mu_X)\big((T')_{\sep (X;\kappa/3,\kappa
  /3)}\big)\geq  m-2\kappa /3$.
  \begin{proof}Suppose that $f_{\ast}(\mu_X)\big((T')_{\sep (X;\kappa/3,\kappa
  /3)}\big)<  m-2\kappa /3$. Then, we have a contradiction since
   \begin{align*}
    \sep \Big(X;\frac{\kappa}{3},\frac{\kappa}{3}\Big)<\dist_T
    \big(T', T\setminus (T')_{\sep (X;\kappa /3,\kappa /3)
    +\varepsilon}\big)\leq \sep
    \Big(f_{\ast}(\mu_X);\frac{\kappa}{3},\frac{\kappa}{3}\Big) \leq
    \sep \Big(X;\frac{\kappa}{3},\frac{\kappa}{3}\Big)
    \end{align*}for any sufficiently small $\varepsilon >0$. 
   \end{proof}
  \end{claim}Combining Claims \ref{cl2} with \ref{cl3}, we obtain
  \begin{align*}
   \diam (f_{\ast}(\mu_X),m-\kappa)\leq  \diam \big((T')_{\sep(X;\kappa
   /3 ,\kappa /3)}\big)
   \leq   2\sep
   \Big(X;\frac{\kappa}{3},\frac{\kappa}{3}\Big)+2\obs_{\mathbb{R}} (X;-\kappa)
   \end{align*}and so (\ref{s16}). This completes the proof of the theorem.
 \end{proof}
Note that the inequality (\ref{s16}) yields slightly worse estimate for the
observable diameter $\obs_{T}(X;-\kappa)$ than (\ref{staba1}). 

Let $T$ be an $\mathbb{R}$-tree and $\nu \in \mathcal{B}^1(T)$ with $m:=\nu(X)$. Taking a median
$m_{\nu}\in T$ of the measure $\nu$, we let $T_{\nu}$ an element in
$\mathcal{C}_{T}'\big(c(\nu)\big)$ with $m_{\nu}\in T_{\nu}$. We then define
  the function $\varphi_{\nu}:T\to \mathbb{R}$ by $\varphi_{\nu} (w):= \dist_T
  (z,w)$ if $w\in T_{\nu}$ and $\varphi_{\nu} (w):= -\dist_T(z,w)$
  otherwise. The function $\varphi_{\nu}$ is clearly the $1$-Lipschitz function.
\begin{lem}\label{p10}Let $T$ be an $\mathbb{R}$-tree and $\nu \in \mathcal{B}^1(T)$. Then,
 the function $\varphi_{\nu}:T \to \mathbb{R}$ satisfies that $c((\varphi_{\nu})_{\ast}(\nu)) \leq 0$,
 \begin{align}\label{s8}
   |c(
  (\varphi_{\nu})_{\ast}(\nu))|\leq \ &
  \crad((\varphi_{\nu})_{\ast}(\nu),m-\kappa)+ \sep
  \Big((\varphi_{\nu})_{\ast}(\nu);\frac{m}{3},\frac{\kappa}{2}\Big)\\
  & \hspace{5.9cm}+\sep
  ((\varphi_{\nu})_{\ast}(\nu);m-\kappa,m-\kappa), \tag*{}
  \end{align}and 
 \begin{align}\label{s9}
  \crad (\nu,m-\kappa)\leq \ &
  \crad((\varphi_{\nu})_{\ast}(\nu),m-\kappa) + \sep \Big(
  \nu;\frac{m}{3},\frac{\kappa}{2}\Big) \\
  & \hspace{1cm}+ \sep
  \Big((\varphi_{\nu})_{\ast}(\nu);\frac{m}{3},\frac{\kappa}{2}\Big)+
  \sep ((\varphi_{\nu})_{\ast}(\nu);m-\kappa,m-\kappa) \tag*{}
  \end{align}for any $\kappa >0$.
 \begin{proof}Combining Lemma \ref{l5} with Proposition \ref{p9}, we have
  \begin{align*}
   \nu(T) c ((\varphi_{\nu})_{\ast}(\nu)) = \int_{T}  \varphi_{\nu} (z) \
   d \nu(z) 
   = \ &\int_{T_{\nu}}\dist_T (c(\nu),z) \ d \nu(z)
   - \int_{T\setminus T_{\nu}}\dist_T (c(\nu),z) \ d \nu (z) \\
   = \ &c_{T_{\nu},c(\nu)} (\nu)\leq 0. 
   \end{align*}Put $r_1:= \crad ((\varphi_{\nu})_{\ast}(\nu),m-\kappa)$ and $r_2:= \sep
  ((\varphi_{\nu})_{\ast}(\nu);m/3,\kappa /2)$. From (\ref{s4}), we observe that $(\varphi_{\nu})_{\ast}(\nu)\big(B_{\mathbb{R}}(\varphi_{\nu}(m_{\nu}),r_2)\big)\geq
  \nu \big( B_T(m_{\nu},r_2)\big)\geq m-\kappa$. Thus, we get
  \begin{align}\label{s10}
   \dist_{\mathbb{R}} \big(  B_{\mathbb{R}}\big( c((\varphi_{\nu})_{\ast}(\nu)),r_1 \big), B_{\mathbb{R}}(
   \varphi(m_{\nu}),r_2)                  \big)\leq \ &\sep ((\varphi_{\nu})_{\ast}(\nu);m-\kappa,m-\kappa)
   \end{align}and so (\ref{s8}). The above inequality
  (\ref{s10}) together with $c((\varphi_{\nu})_{\ast}(\nu))\leq
  0$ yields that
  \begin{align*}
   \dist_T (c(\nu),m_{\nu})= \varphi_{\nu} (m_{\nu})\leq \ &|c(
   (\varphi_{\nu})_{\ast}(\nu))-\varphi_{\nu}(m_{\nu})|\\
   \leq \ &r_1 +r_2 + \sep ((\varphi_{\nu})_{\ast}(\nu);m-\kappa,m-\kappa)=:r_3.
   \end{align*}Therefore, putting $r_4:= \sep (\nu ; m/3, \kappa /2)$, we obtain
  \begin{align*}
   \nu \big(B_T(c(\nu),r_3+r_4)\big) \geq
   \nu \big(B_T(m_{\nu},r_4)\big)\geq m-\kappa
   \end{align*}and so (\ref{s9}). This completes the proof.
  \end{proof}
 \end{lem}

\begin{prop}\label{c3}Let $T$ be an $\mathbb{R}$-tree and $X$ an mm-space with
 $\mu_X \in \mathcal{B}^1(X)$. Then, for any $\kappa >0$ we have
 \begin{align*}
  \obsc_T(X;-\kappa) \leq \obsc_{\mathbb{R}}(X;-\kappa)+2\sep
  \Big(X;\frac{m}{3}, \frac{\kappa}{2}\Big) + \sep (X;m-\kappa,m-\kappa).
  \end{align*}
 \begin{proof}
  This follows from Lemma \ref{l1} and Lemma \ref{p10}.
  \end{proof}
 \end{prop}

 \begin{proof}[Proof of Theorem \ref{t1}]Proposition \ref{c3} together
  with Corollary \ref{c1} and Lemma \ref{l7} directly implies the proof of the theorem. 
  \end{proof}

  \begin{lem}\label{p11}Let $T$ be an $\mathbb{R}$-tree and $\nu \in
   \mathcal{B}^1(T)$. Then, for any $p\geq 1$ and $\kappa >0$, we have
   \begin{align}\label{sup1}
    V_p(\nu) \leq \ & 2m^{2/p} \Big\{ \crad(
    (\varphi_{\nu})_{\ast}(\nu),m-\kappa) + \sep \Big(
    (\varphi_{\nu})_{\ast}(\nu);\frac{m}{3},\frac{\kappa}{2}\Big) \\
    & \hspace{4.5cm}      + \sep (
    (\varphi_{\nu})_{\ast}(\nu);m-\kappa,m-\kappa) \Big\} +2V_p(\varphi_{\nu}).
    \tag*{}
   \end{align}In the case of $p=2$, we also have the better etimate
    \begin{align}\label{sup2}V_2(\nu )^2 \leq \ &4m^2 \Big\{\crad ( (\varphi_{\nu})_{\ast}(\nu),m-\kappa) +   \sep
  \Big((\varphi_{\nu})_{\ast}(\nu);\frac{m}{3},\frac{\kappa}{2}\Big)\\ \
    & \hspace{4cm} +\sep
  ((\varphi_{\nu})_{\ast}(\nu);m-\kappa,m-\kappa)      \Big\}^2
     + 2V_2(\varphi_{\nu})^2. \tag*{}
    \end{align}
   \begin{proof}From the triangle inequality, we have
    \begin{align}\label{s11}
     V_p(\nu) \leq 2 \Big(   \int \int_{T\times T} \dist_T
     (c(\nu),z)^p       \ d\nu (z)d\nu (w)  \Big)^{1/p}  = 2 \Big(m
     \int_T \dist_T (c(\nu),z)^p       \ d\nu (z)\Big)^{1/p}.
     \end{align}Putting
    $c_{\nu}:=c((\varphi_{\nu})_{\ast}(\nu))$, we also get
    \begin{align}\label{s12}
     \Big(\int_T \dist_T (c(\nu),z)^p       \ d\nu (z) \Big)^{1/p}
     = \ & \Big(\int_T |\varphi_{\nu}(z)|^p  \
     d\nu(z) \Big)^{1/p}  \\
\leq \ &   m^{1/p} |c_{\nu}|  +  \Big(\int_{\mathbb{R}} |c_{\nu}-r|^p \
     d(\varphi_{\nu})_{\ast}(\nu)(r) \Big)^{1/p}  \tag*{} \\
     \leq \ & m^{1/p}|c_{\nu}| + \frac{V_p(\varphi_{\nu})}{m^{1/p}}, \tag*{}
     \end{align}where in the last inequality we used Lemma \ref{l5}. 
    Combining (\ref{s11}) with (\ref{s12}), we obtain (\ref{sup1}).

    In the case of $p=2$, we have
    \begin{align}\label{sd1}
   \int_T \dist_T(c(\nu),z)^2
 \ d\nu(z) =\ & \int_{\mathbb{R}} |r|^2 \
     d(\varphi_{\nu})_{\ast}(\nu)(r) \\
     = \ & m|c_{\nu}|^2 + \int_{\mathbb{R}} |r-c_{\nu}|^2 \
     d(\varphi_{\nu})_{\ast}(\nu)(r) \tag*{}\\
     = \ & m|c_{\nu}|^2 + \frac{V_2(\varphi_{\nu})^2}{2m}, \tag*{}
    \end{align}where in the second and the last equalities we used Lemma
    \ref{l5}. Substituting (\ref{sd1}) to (\ref{s11}), we obtain
    (\ref{sup2}). This completes the proof.
   \end{proof}
   \end{lem}

   \begin{prop}\label{c4}Let $T$ be an $\mathbb{R}$-tree and $X$ an
    mm-space. Then, for any $p\geq 1$, we have
    \begin{align}\label{sup5}
     \obin_T(X) \leq 2 \big\{2^{1/p}(1+ 2\cdot 2^{1/p})+1 \} \obin_{\mathbb{R}}(X). 
     \end{align}In the case of $p=2$, we also have the better estimate
    \begin{align}\label{s15}
     \obinin_{T}(X)^2 \leq (38+16\sqrt{2})\obinin_{\mathbb{R}}(X)^2. 
     \end{align}
    \begin{proof}Assume first that $f_{\ast}(\mu_X)\in
     \mathcal{B}^1(T)$ for any $1$-Lipschitz map $f:X\to T$. Then,
      Lemma \ref{l1} together with Lemma \ref{l2} and (\ref{sup1}) implies
     that
     \begin{align*}
      \obin_T(X)\leq \ &2m^{2/p} \Big\{  \obsc_{\mathbb{R}} (X;-\kappa) +\sep \Big(
      X;\frac{m}{3},\frac{\kappa}{2}\Big)\Big\} +
      2\obin_{\mathbb{R}}(X) \\
      \leq \ & 2m^{2/p} \Big\{ \obsc_{\mathbb{R}}(X;-\kappa) + \sep
      \Big( X;\frac{\kappa}{2}, \frac{\kappa}{2}\Big)\Big\} + 2\obin_{\mathbb{R}}(X)
      \end{align*}for any $0<\kappa < m/2$. Hence, applying the
     inequalities (\ref{exc3}) and (\ref{exc5}) to this inequality, we
     get
     \begin{align*}
      \obin_T(X)\leq 2\big\{ m^{1/p}\kappa^{-1/p}( 1+2\cdot 2^{1/p}   )
      + 1 \big\}\obin_{\mathbb{R}}(X)
      \end{align*}for any $0<\kappa < m/2$. Letting $\kappa \to m/2$, we
     get (\ref{sup5}). In the case of $p=2$, from (\ref{sup2}), we have
     \begin{align*}
      \obinin_{T}(X)^2 \leq 4m^2 \Big\{ \obsc_{\mathbb{R}}(X;-\kappa) +
      \sep \Big(X;\frac{\kappa}{2},\frac{\kappa}{2}\Big) \Big\}^2 + 2\obinin_{\mathbb{R}}(X)^2
      \end{align*}for any $0<\kappa<m/2$. Therefore, substituting the
     inequalities (\ref{exc4}) and (\ref{exc6}) to this inequality, we get
     \begin{align*}
      \obinin_{T}(X)^2 \leq 2\big\{m\kappa^{-1}  (2\sqrt{2}+1)^2         +1\big\}  \obinin_{\mathbb{R}}(X)^2
      \end{align*}for any $0<\kappa <m/2$. Letting $\kappa \to m/2$, we obtain (\ref{s15}).

     We consider the other case that there exists a $1$-Lipschitz map
     $f:X\to T$ with $f_{\ast}(\mu_X)\notin \mathcal{B}^1(T)$. By using
     H\"{o}lder's inequality and Fubini's theorem, we have $V_p(f)=+\infty$. Taking $x_0 \in X$,
     we put $f_n:= f|_{B_X(x_0,n)}$ for each $n\in \mathbb{N}$. From
     Lemma \ref{l3} and the
     above proof, we have
     \begin{align*}
      V_p(f_n) \leq \obin_{T}\big(B_X(x_0,n)\big) \leq \ &
      2\{ 2^{1/p}(1+ 2\cdot 2^{1/p}) +1    \}\obin_{\mathbb{R}}\big(B_X(x_0,n)\big) \\
      \leq \ &2\{ 2^{1/p}(1+ 2\cdot 2^{1/p}) +1    \}\obin_{\mathbb{R}}(X).
      \end{align*}Since $V_2(f_n)\to V_2(f)=+\infty$ as $n\to \infty$, this
     implies $\obin_{\mathbb{R}}(X)=+\infty$. This completes the proof. 
     \end{proof}
    \end{prop}

    \begin{proof}[Proof of Theorem \ref{t2}]Proposition \ref{c4} directly implies the proof
     of the theorem..
     \end{proof}

 \begin{ack}\upshape The author would like to express his thanks to
  Professor Takashi Shioya for his valuable suggestions and assistances during the preparation of
  this paper. He also thanks Professor Vitali Milman and Professor
  Shin-ichi Ohta for useful comments.
  \end{ack}

	\end{document}